\documentclass[11pt]{article}
\usepackage{authblk} 

\usepackage[latin1]{inputenc}
\usepackage[T1]{fontenc}
\usepackage{amsmath}
\usepackage{amsfonts}
\usepackage{amssymb}
\usepackage{enumerate}
\usepackage[normalem]{ulem} 

\usepackage{amsmath,xcolor,ulem}
\usepackage{amsfonts}
\usepackage{amssymb,todonotes}
\usepackage{amsthm}
\usepackage{graphicx}
\usepackage{epsfig}
\usepackage{subcaption}
\usepackage{marginnote}

\usepackage[colorlinks=true, allcolors=blue]{hyperref}

\usepackage[top=2.8cm,bottom=2.8cm,left=2.5cm,right=2.5cm]{geometry}
\usepackage{float}
\usepackage[algoruled]{algorithm2e}
\usepackage[font=footnotesize,width=\linewidth]{caption}

\newcommand{\F}{\mathcal{F}}
\newcommand{\Dtau}{\Delta \tau}

\newcommand{\tY}{\Tilde{Y}}
\newcommand{\half}{\frac{1}{2}}

\newcommand*{\bm}[2]{W_{#1}^{#2}}

\newcommand{\Vd}[0]{V_\text{data}}

\newcommand{\sn}[0]{\sigma_\nu}
\newcommand{\kn}[0]{\kappa_\nu}
\newcommand{\mr}[0]{\mu_{\nu}}

\newcommand{\tOmega}[0]{\tilde{\Omega}}
\newcommand{\mean}[0]{\text{mean}}

\DeclareMathOperator{\argmin}{argmin}
\floatname{algorithmic}{Procedure}

\title{A Space Mapping approach for the calibration of financial models with the application to the Heston model}
\author{Anna Clevenhaus$^1$, Claudia Totzeck$^2$ and Matthias Ehrhardt$^1$}
\affil{$^1$Applied and Computational Mathematics, University of Wuppertal, Germany}
\affil{$^2$Continuous Optimization, University of Wuppertal, Germany}

\begin{document}
\maketitle

\begin{tikzpicture}[remember picture,overlay]
	\node[anchor=north east,inner sep=20pt] at (current page.north east)
	{\includegraphics[scale=0.2]{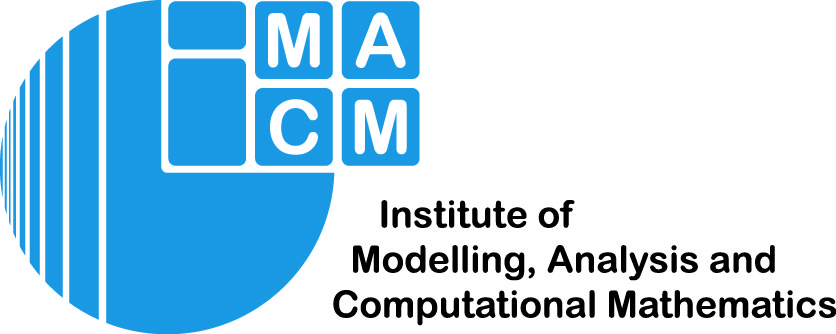}};
\end{tikzpicture}

\begin{abstract}
We present a novel approach for parameter calibration of the Heston model for pricing an Asian put option, namely space mapping. 
Since few parameters of the Heston model can be directly extracted from real market data, calibration to real market data is implicit and therefore a challenging task. 
In addition, some of the parameters in the model are non-linear, which makes it difficult to find the global minimum of the optimization problem within the calibration. 
Our approach is based on the idea of space mapping, exploiting the residuum of a coarse surrogate model that allows optimization and a fine model that needs to be calibrated.
In our case, the pricing of an Asian option using the Heston model SDE is the fine model, and the surrogate is chosen to be the Heston model PDE pricing a European option. 
We formally derive a gradient descent algorithm for the PDE constrained calibration model using well-known techniques from optimization with PDEs.
Our main goal is to provide evidence that the space mapping approach can be useful in financial calibration tasks. 
Numerical results underline the feasibility of our approach.
\end{abstract}

\begin{minipage}{0.9\linewidth}

\medskip

\end{minipage}

%%%%%%%%%%%%%%%%%%%%%%%%%%%%%%%%%%%%%%%%%%%%%%%%%%%%%%%%%%%%%%%%%%%%% INTRO
\section{Introduction}
Financial models aim to approximate the real market behavior of various underlyings. 
Since most model parameters are implicit in the real market data, calibrating the model parameters to fit the real market data is challenging \cite{Cui.2017,horvath2021deep,liu2019neural,mehrdoust2023two,Mikhailov.2003,Teng.2016}. 
In this article, we consider the well-known Heston model \cite{Heston.1993}, a two-dimensional stochastic differential equation (SDE) that allows to simulate the behavior of stock prices. 
In particular, for a given variance of the option we calibrate four implicit parameters -- the volatility of the variance, the mean inversion rate, the long-term mean and the correlation factor of the Brownian motions in the model equations for the asset and the variance -- to the market data. 
Our goal is to provide a proof-of-concept, so we choose the simplest model with constant coefficients.
However, we are aware that the Heston model can be extended in various ways, e.g.\ by considering time-dependent parameters \cite{Mikhailov.2003} or by including additional processes, which leads to an increase in the difficulty of the calibration process \cite{Teng.2016}.

Various calibration techniques for the Heston model are available in the literature. 
Indeed, in the setting of constant parameters and for very specific use cases, there are approaches based on the closed-form valuation formula of the Heston model \cite{Mikhailov.2003,Cui.2017}. 
These are fast and provide information about the global minimum. 
For more general cases, the stochastic nature of the Heston model allows for Monte Carlo optimization methods \cite{Teng.2016}, which can also be used to calibrate the stock price and variance. 
The Monte Carlo theory is well established, but the approaches are computationally expensive and do not provide information about global or local minima. 
Recently, calibration approaches using neural networks, deep learning strategies as well as parallel GPU implementation of the Heston model have been proposed \cite{horvath2021deep, liu2019neural, Liu.2019, Leite.2021, Ferreiro.2020}.
The networks must be trained individually for each model, and training requires appropriate data. Again, there is no information about global or local minima. 
An advantage is that once the neural networks are trained, they can be evaluated quickly.

%%%%%%%%%%%%
Our contribution to this line of research is a calibration algorithm for the Heston model that is independent of a specific characteristic function and easily extendable to time-dependent parameters.
The core of the algorithm is based on the space mapping \cite{bandler2004space}, an iterative procedure that minimizes the residuum of a fine and a coarse model.
In fact, to calibrate the parameters of the fine model, the coarse model is optimized and the fine model is only evaluated. 
To demonstrate the proposed approach, we calibrate the Heston SDE model for pricing Asian put options (fine model) to real market data. The calibration uses only one given market price. 
In our setting, the coarse model is chosen because it is the computationally less expensive Heston model for pricing the European put option, which is known to be a good approximation for the Asian option in the case of short maturities. 
This is because the European option is a limiting case for the Asian option. 
In more detail, while we aim to calibrate the parameters of the SDE for the Asian put option, we will optimize the deterministic PDE model of the European option, which can be solved using techniques from optimization with PDE \cite{HUUP,Troeltzsch}. 
In each iteration of the calibration process, we will evaluate the SDE to compute the residuum of the two models. 
To our knowledge, there is only one paper dealing with the calibration of Asian options under the Heston model, where Khalife et al.~\cite{Khalife.2023} calibrate the Heston model using artificial intelligence and the approximation is done via the European option as well.
%In our work we use the space mapping approach to calibrate our Heston model using only one market data value.
%We consider the Heston model in its $\log$-transformed form, as the control model for the space mapping approach. 
%The space mapping algorithm uses a coarse model to approximate the fine model calibration to solve the control problem \cite{bandler2004space}.
%The fine model is given by the stochastic differential equation (SDE) representation of the Heston model for pricing the Asian put option, which is approximated by the corresponding partial differential equations (PDEs) for the Heston model for pricing the European put option, as the coarse model.

%Within this space mapping approach, we have two differences between the coarse and fine models. 
%One difference is the choice of option type. 
%On the one hand, the choice is motivated by the fact that the computational effort for pricing European options is less than for Asian options. 
%This is because European options only consider the price of the underlying asset at expiration, instead of considering the value at each time step, as Asian options do. 
%This approximation is also used in \cite{Khalife.2023}.
%The second difference is the representation of the model, as it is more difficult to calibrate an SDE model instead of a deterministic PDE. 
%For the space mapping approach, we consider the \textit{aggressive space mapping} (ASM) approach . 

%The Heston model contains at least four parameters implicit in the market data, in our case we consider a calibration of five parameters. 
The proposed calibration is based on \textit{aggressive space mapping} (ASM) \cite{Echeverria.2008}.
 To iteratively update the parameters of the fine model, a PDE-constrained optimization problem is considered at the level of the coarse model. 
We measure the difference between the fair price given by the numerical solution of our model and the reference data, the subsequent market data, in the cost functional. 
The PDE-constrained optimization problem is solved using a gradient descent method. 
For this purpose, we formally derive an adjoint-based gradient descent algorithm for the Heston PDE model.
 Note that gradient descent algorithms have previously been used implicitly in the context of neural network approximations for the Heston model \cite{Liu.2019}.

As mentioned above, our main goal is to provide a proof-of-concept that the space mapping approach can be useful for calibration tasks in financial research.
 This claim is supported by the numerical results. 
However, we are aware that there is much room for improvement, for example an extension to time-dependent parameters, see \cite{clevenhaus24gradient} for details on the PDE-constrained optimization problem. 

The article has the following structure.
In Section~\ref{sec: Financial modelling}, we present the mathematical formulation of the Heston SDE model and its $\log$-transformed form for pricing an Asian put option.
We also discuss the PDE formulation of the Heston model for pricing a European put option.
In Section~\ref{sec:spacemapping}, we then introduce the corresponding optimization problems and present details of the  space mapping approach. This includes the formal derivation of the adjoint of the $\log$-transformed Heston model as well as the identification of the gradient which is used in the descent algorithm. 
We then discretize the SDE model and the PDE calibration method in Section~\ref{sec: Numerical scheme}. 
Numerical results to illustrate the feasibility of space mapping in finance are presented in Section~\ref{sec: Numerical Results}.
Finally, in Section~\ref{sec: Conclusion and Outlook} the paper ends with a conclusion and an outlook.

%%%%%%%%%%%%%%%%%%%%%%%%%%%%%%%%%%%%%%%%%%%%%%%%%%%%%%%%% Section 2
\section{The Heston SDE model and its PDE approximation}\label{sec: Financial modelling}
The Heston model was proposed by Heston in 1993 \cite{Heston.1993} and describes the dynamics of the underlying asset and its variance by a two-dimensional SDE. 
It is an extension of the well-known Black-Scholes model \cite{Black.1973}. 
Let $S\in \mathbb{R}^{+}$ denote the asset, $r\in \mathbb{R}^{>0}$ the risk-free interest rate, $\sigma_S\in \mathbb{R}^{>0}$ the volatility and $\bm{t}{S}$ the Brownian motion.  
Then the Black-Scholes model describes the dynamics of $S$ as a geometric Brownian motion process given by
\begin{equation}
    dS_t = r S_t\,dt + \sigma_S S_t\,d\bm{t}{S},  \quad t \in [0,T],\; S_0>0,\\
\end{equation}
where $T$ is the maturity.

Heston introduced a second stochastic process to model more complex volatility behavior.
 Let $\nu\in \mathbb{R}^{+}$ denote the variance, set $\sigma_S=\sqrt{\nu}$, and assuming that the volatility follows an Ornstein-Uhlenbeck process, he includes a Cox-Ingersoll (Ross) process for the variance. 
 This leads to Heston's SDE model under the risk-neutral measure given by
\begin{equation}\label{eq:pureHeston}
 \begin{cases}
 dS_t = (r-q) S_t\,dt + \sqrt{\nu_t} S_t\,d\bm{t}{S},&\quad S_0>0,\\
 d\nu_t = \kn(\mr-\nu_t)\,dt+\sn \sqrt{\nu_t}\,d\bm{t}{\nu},&\quad \nu_0>0,
\end{cases},
\end{equation}
where $q\in\mathbb{R}^{>0}$ is the continuous dividend rate, $\kn\in\mathbb{R}^{>0}$ is the mean reversion rate, $\mr\in \mathbb{R}^{>0}$ is the long term mean,
and $\sn\in\mathbb{R}^{>0}$ is the volatility of the variance. 
The Brownian motions $\bm{t}{S}$ and $\bm{t}{\nu}$ are correlated by the constant $\rho \in [-1,1]$, cf.\ \cite{Heston.1993}.
To ensure that the square root of \eqref{eq:pureHeston} is positive, the Feller condition $2\kn \mr \ge\sn^2$ must be satisfied. 
Otherwise, computational problems arise due to a complex variance.
The Heston SDE representation is used to price Asian put options with the payoff function
\begin{equation}
      \phi_S(S,T)=\max \bigl( K -  A_S(0,T)\bigr),
\end{equation}
where $K$ is the predetermined strike price and
\begin{equation}
    A_S(0,T)=\frac{1}{T} \int_{0}^{T} S_t \,dt.
\end{equation}

%%%%%%%%%%%%%%%%%%%%%%%%%%%%%%%%%%%%%%%%%%%%%%%%%%%%%%
\subsection{The fine model}
To zoom into the particularly interesting price range near $K$, we use the variable transformation $x=\log(S)$ for the asset. We get the \textit{$\log$-transformed Heston SDE}
\begin{equation}\label{eq:LogFine}
\begin{cases}
  dx_t = (r- q - \half \nu_t)\,dt + \sqrt{\nu_t} \,d\bm{t}{S},&\quad S_0>0,\\
  d\nu_t = \kn(\mr-\nu_t)\,dt+\sn \sqrt{\nu_t}\,d\bm{t}{\nu},&\quad \nu_0>0, \\
  d\bm{t}{x}\,d\bm{t}{\nu}=\rho\,dt
\end{cases}.
\end{equation}
on the semi-unbounded domain $x\in\mathbb{R}$, $\nu\ge0$, $0\le t \le T$ with the transformed payoff function
\begin{equation}\label{eq: payoff x}
    \phi(x,T)=\max \bigl( K - A(0,T) \bigr)
    = \max \Bigl( K- \frac{1}{T} \int_{0}^{T} \exp(x_t) \,dt\Bigr).
\end{equation}
The payoff condition holds at the maturity and we are interested in the fair price $V_f(x,t)$ today at $t=0$ given by
\begin{equation}\label{eq: V_f}
    V_f(x,0)=\exp(-rT) \,\mathbb{E}\bigl(\phi(x,T)\bigr)
\end{equation}
using the discounted expectation value. 
This model will be used as the fine model in the space mapping approach discussed in Section~\ref{sec:spacemapping}. 

%%%%%%%%%%%%%%%%%%%%%%%%%%%%%%%%%%%%%%%%%%%%%%%%%
\subsection{The coarse model}
As we will discuss in more detail below, the theory of space mapping requires that the coarse model be a reasonable approximation of the fine model that can be easily optimized.
Since the European plain vanilla option is a limiting case for the Asian put option, we choose it as the coarse model.
Furthermore, we use the PDE formulation of the Heston model for the European plain vanilla option because we can optimize it using methods from optimization with PDE constraints. 

Note that the payoff condition \eqref{eq: payoff x} holds at maturity and thus yields a terminal condition. 
To obtain an initial value problem, we reverse the time $\tau=T-t$.  
Then, the Heston PDE under a risk-neutral measure is derived using Kolmogorov's backward equation. 
For the fair price of an option $V(x,\nu,\tau)$ it reads
\begin{equation}\label{eq: Heston transformed}
\begin{split}
    V_{\tau}= \frac{\nu}{2} V_{xx}+ \frac{1}{2} \sigma_{\nu}^2 \nu V_{\nu \nu}+ (r-q-\frac{\nu}{2})V_x+\kn (\mr-\nu)V_{\nu}+ \sigma_{\nu} \nu \rho V_{x \nu} -rV ,
\end{split}
\end{equation}
and is supplemented with the initial condition 
\begin{equation}
   V(x,\nu,\tau)=\max\bigl(K-\exp(x),0\bigr),
   % V(x,\nu,\tau)=\bigl(K-\exp(x),0\bigr)^+
\end{equation}
and appropriate boundary conditions. 
We choose to use the boundary conditions proposed by Heston and apply the $\log$-transformation to obtain
\begin{align}
    V(x,\nu,\tau)&\sim K\exp(-r\tau), \quad\text{for }x\to-\infty,\\
    V(x,\nu,\tau)&=0,\quad\text{for }x\to\infty,\\
    V(x,\nu,\tau)&\sim K\exp(-r\tau),\quad\text{for }\nu\to\infty.\label{bc_nu_infty}
\end{align}

Note that the parabolic PDE \eqref{eq: Heston transformed} degenerates to a first-order hyperbolic PDE at $\nu=0$, therefore we need to consider the Fichera theory \cite{Buckova.2016, Kutik.2015} to assess whether it is necessary to provide an analytic boundary condition at $\nu=0$ or not.  
From the Fichera condition at $\nu=0$ given by 
\begin{equation}\label{eq:fichera}
   b(\nu)= \kn ( \mr-\nu) - \half \sn^2 
\end{equation}
we get a case dependency. 
At $\nu=0$ we get an outflow boundary if $\lim_{\nu \to 0^+} b(\nu)\ge0$, otherwise we get an inflow boundary. 
The condition for an outflow boundary \eqref{eq:fichera} is equivalent to the Feller condition $2\kn\mr\ge\sn^2$, which motivates us to assume in the following that the Feller condition holds, so that we do not have to supply an analytical boundary condition at $\nu=0$. At this point we want to mention that we need a numerical boundary condition to complete the scheme, which will be discussed later in Section~\ref{sec: Gradient Descent Algorithm}.

%%%%%%%%%%%%%%%%%%%%%%%%%%%%%%%%%%%%%%%%%%%%%%% Section 3
\section{The Space Mapping Approach}\label{sec:spacemapping}
Below we summarize the main ideas of the space mapping approach. For more details we recommend \cite{bandler2004space,Echeverria.2008,Totzeck2019SpaceMR,Marheineke2012SpaceMC}. 
The purpose of space mapping is to optimize (or calibrate) an accurate and computationally expensive model (fine) that cannot be optimized using a surrogate model (coarse) for which efficient optimization algorithms are available.

In our case, the accurate model is the Heston SDE for Asian option pricing, and the coarse model is the Heston PDE for European option pricing, as discussed above.
%In this way, we obtain a surrogate model structure and reduce the computational time in optimization, where time-consuming models are needed to obtain sufficiently accurate results.
Both models return a single price for a predefined asset $\tilde{S}$ and a variance $\tilde{\nu}$.  
The predefined asset is part of the contract, but $\tilde{\nu}$ must be determined implicitly from the market data or, as in our case, predefined by a guess. 
Together with the other unknown model parameters, we want to calibrate $\xi=(\sn,\rho,\kn,\mr)$. 
Note that we can take advantage of the fact that the parameters we want to calibrate are present in both models. 

Although we expect the optimal values of the parameters to be different for the two models, the space mapping technique helps us to calibrate the parameters of the fine model while only optimizing the coarse model. 
In the following, we distinguish between $\xi_{\mathrm{f}}$, the parameter vector for the fine model, and $\xi_{\mathrm{c}}$, the parameter vector for the coarse model. Similarly, we denote the option price of the fine model by $V_f$ and the option price of the coarse model by $V_c$.

Since we are using real market data $\Vd$ as ground truth for the calibration, we define the cost functional as follows
\begin{equation}\label{eq: J}
   J(V;\Vd) = \frac{1}{2} \int_0^T \| V - \Vd \|^2 \,d\tau.
\end{equation}

In the following, we want to approximate the solution of the fine calibration problem
\begin{align}\label{prob:fine}
\min_{\xi_{\mathrm{f}} \in X_{\mathrm{f}}} J(V_f;\Vd) \quad
\text{subject to} \quad V_f(x,0) = \exp(-rT)\mathbb{E}(\phi(x,T)), 
\end{align}
where $x$ is the solution of the log-transformed Heston SDE \eqref{eq:LogFine}.
 Note that we do not solve the fine optimization problem, but only evaluate the cost functional for given parameter sets $\xi_f$ during the space mapping algorithm.

%The fine model is given by $f(\xi_{\mathrm{f}})\in \mathbb{R}$, where $\xi_{\mathrm{f}} \in X_{\mathrm{f}} \in \mathbb{R}^{5}$ is the fine model variable, as in our case $f(\xi_{\mathrm{f}})$ is the option price. 
%Therefore, we can define a \textit{fine model cost functional} 
%\begin{equation}
%     \min_{\xi_{\mathrm{f}} \in X_{\mathrm{f}}} J(V_f,f(\xi_{\mathrm{f}})) 
%\end{equation}
%where 
%\begin{equation}\label{eq: J_f}
%   J_f(V_f, f(\xi_{\mathrm{f}})) = \frac{1}{2} \int_0^T \| V_f - \Vd \|^2 \,d\tau.
%\end{equation}
%measures the difference between the market data $V$ and the particular response of the mathematical problem. The goal is to minimize the cost functional
%\begin{equation}
%    \xi_{\mathrm{f}}^*:= \argmin_{\xi_{\mathrm{f}} \in X_{\mathrm{f}}} J_f(V_f,f(\xi_{\mathrm{f}})).
%\end{equation}

The optimization problem, which aims to calibrate the coarse model for given real market data, is given by 
\begin{equation}\label{prob:coarse}
\min_{\xi_{\mathrm{c}} \in X_{\mathrm{C}}} J(V_c;\Vd) \quad
\text{subject to} \quad V_c \text{ solves } \eqref{eq: Heston transformed}.
\end{equation}
In the following we assume that both problems \eqref{prob:fine} and \eqref{prob:coarse} admit a unique solution with a unique minimizer.

The solution of the coarse calibration problem will be our initial guess for the optimal parameter set of the fine model.
However, in the space mapping approach we want to iteratively improve the parameter values of the fine model and we exploit the approximation properties of the coarse model. Let's assume that there exists a so-called \textit{space mapping function} $s\colon X_{\mathrm{f}}\to X_{\mathrm{c}}$, $\xi_{\mathrm{f}}\mapsto s(\xi_{\mathrm{f}})$, which satisfies
\begin{equation}\label{eq:S}%\tag{S}
s(\xi_{\mathrm{f}}):= \argmin_{\xi_{\mathrm{c}}\in X_{\mathrm{c}}} r\bigl(V_c(\xi_{\mathrm{c}}), V_f(\xi_{\mathrm{f}})\bigr),
\end{equation}
for some misalignment function $r\colon Y\times Y\to \mathbb{R}$. %In the following we choose $r(V_c,V_f) = \frac12 \int_0^T \| V_c - V_f \|^2 d\tau$.

Assuming that the coarse model is a good approximation of the fine model, and choosing $J$ as the misalignment function, we expect the following condition to be satisfied
\begin{align*} %\label{eq:S2}\tag{S2}
  s(\xi_{\mathrm{f}}^*) &= \argmin_{\xi_{\mathrm{c}}\in X_{\mathrm{c}}} J\bigl(V_c(\xi_{\mathrm{c}}); V_f(\xi_{\mathrm{f}}^*)\bigr) \\
  &\approx
  \argmin_{\xi_{\mathrm{c}}\in X_{\mathrm{c}}} J\bigl(V_c(\xi_{\mathrm{c}}); V_c(\xi_{\mathrm{c}}^*)\bigr) \\
  &\approx
  \argmin_{\xi_{\mathrm{c}}\in X_{\mathrm{c}}} J\bigl(V_c(\xi_{\mathrm{c}}); \Vd\bigr)
  = \xi_{\mathrm{c}}^*.
\end{align*}
 The underlying assumption is, that the optimal states $V_f$ and $V_c$ are both good approximations of the ground truth $\Vd$, each for the respective model.

The strategy of the space mapping algorithm is to solve
\begin{equation*}% \label{eq:S3}\tag{S3}
    s(\xi_{\mathrm{f}}^*)- \xi_{\mathrm{c}}^*=0.
\end{equation*}
Note that we will not approximate the entire function $s$, but only evaluate it along the iterations.
Approximating $s$ entirely is a much harder and probably ill-posed task.

For the numerical results we use a simplified version of the Aggressive Space Mapping (ASM) algorithm \cite{Marheineke2012SpaceMC}, in fact by the linearity of the state problems \eqref{eq:LogFine} and \eqref{eq: Heston transformed} we can approximate the Jacobian of the space mapping function by the identity \cite{Bark.2000} and thus obtain the Algorithm~\ref{alg:SpaceMappingsimple}.

% \bigskip
% \begin{algorithm}[H]\label{alg:SpaceMapping}
% 	\SetAlgoLined
% 	\KwResult{optimized $\xi_{\mathrm{f}}$}
% 	$\xi_{\mathrm{f}}^0=\xi_{\mathrm{c}}^*=\argmin_{\xi_{\mathrm{c}}\in X_{\mathrm{c}}} \| c(\xi_{\mathrm{c}}) - V \|$\;
%     $B_0= I$\;
%     \For{$k=0,1,\ldots$}{
%         1) evaluate space mapping function $s(\xi_{\mathrm{f}}^k)$ by \;
%          $\quad$ (a) evaluate fine model \eqref{eq: V_f} to obtain $V_f$ using $\xi_{\mathrm{f}}^k$ \\
%          $\quad$ (b) Perform a coarse model optimization \\$\qquad s(\xi_{\mathrm{f}}^k)=\argmin_{\xi_{\mathrm{c}}\in X_{\mathrm{c}}} \| c(\xi_{\mathrm{c}}) - V_f \|$\\
%          2) Evaluate $B^k=B^{k-1}+ ((s(\xi_{\mathrm{f}}^k)-\xi_{\mathrm{c}}^*) \otimes h^{k-1} )/(h^{k-1})^2$ for $k>0$ \\
%          3) Solve $B^k h^k= -(s(\xi_{\mathrm{f}}^k)-\xi_{\mathrm{c}}^*))$ for update $h_k$\\
%          4) Update control $\xi_f^{k+1}=\xi_f^k+h^k$ using the projected Armijo rule to restrict $xi_f^{k+1}$ to the boundary\\
%          while 	$\| s(\xi_{\mathrm{f}}^k) - \xi_{\mathrm{c}}^*\| >$ tolerance\;
%     }
% 	\caption{The Aggressive Space Mapping (ASM) algorithm.}
% \end{algorithm}
% \bigskip
% \noindent
% where the space mapping function is approximated by linearization using a Broyden-type approximation for the Jacobian of the space mapping function, see Algorithm~\ref{alg:SpaceMapping}.
% Since we can assume that the mapping is linear, we can set $B=I$ as a constant \cite{Bark.2000}.
% Then the ASM algorithm simplifies to Algorithm~\ref{alg:SpaceMappingsimple}.

\bigskip
\begin{algorithm}[H]\label{alg:SpaceMappingsimple}
	\SetAlgoLined
	\KwResult{optimized $\xi_{\mathrm{f}}$}
	$\xi_{\mathrm{f}}^0=\xi_{\mathrm{c}}^*=\argmin_{\xi_{\mathrm{c}}\in X_{\mathrm{c}}} \| c(\xi_{\mathrm{c}}) - V \|$\;
    \For{$k=0,1,\ldots$}{
        1) evaluate space mapping function $s(\xi_{\mathrm{f}}^k)$ by \;
         $\quad$ (a) evaluate fine model \eqref{eq: V_f} to obtain $V_f$ using $\xi_{\mathrm{f}}^k$ \\
         $\quad$ (b) Perform a coarse model optimization \\$\qquad s(\xi_{\mathrm{f}}^k)=\argmin_{\xi_{\mathrm{c}}\in X_{\mathrm{c}}} J\bigl(V_c(\xi_{\mathrm{c}}); V_f(\xi_{\mathrm{f}})\bigr)$\\
         2) Compute $h^k= -(s(\xi_{\mathrm{f}}^k)-\xi_{\mathrm{c}}^*))$ for update $h_k$\\
         4) Update control $\xi_f^{k+1}=\xi_f^k+h^k$ using the projected Armijo rule to restrict $\xi_f^{k+1}$ to the boundary
         \\while
    	$\| s(\xi_{\mathrm{f}}^k) - \xi_{\mathrm{c}}^*\| >$ tolerance\;
    }
	\caption{The simplified Aggressive Space Mapping (ASM) algorithm}
\end{algorithm}
\bigskip

%%%%%%%%%%%%%%%%%%%%%%%%%%%%%%%%%%%%%%%%%%%%%%% Section 3.1
\subsection{Optimal Control of the Coarse Model}\label{sec: Gradient Descent Algorithm}
In the following, we formally derive the first-order optimality system of the coarse calibration problem that we solve in (b). 
For notational convenience, we use $J\bigl(V_c(\xi_{\mathrm{c}}); \Vd)\bigr)$ in the derivation. Our formal derivation is based on a Lagrangian.
Specifically, we formally derive a gradient-based algorithm using a Lagrangian approach to solve  
\begin{equation}
     \xi_{\mathrm{c}}^*:= \argmin_{\xi_{\mathrm{c}} \in X_{\mathrm{c}}} J(V_c,\Vd) \text{ subject to } \eqref{eq: Heston transformed}
\end{equation}
for a given data set $\Vd$. 
We denote the Lagrange multipliers by $\psi = (\varphi, \varphi^a, \varphi^b, \varphi^c,\varphi^d)$, set $\Omega = (0,\nu_{\max}) \times (x_{\min}, x_{\max})$ and partition the boundary $\partial \Omega$ into 
\begin{align*}
 &\Gamma_a = \partial\Omega \cap \{x=-\infty\}, && \Gamma_b = \partial\Omega \cap \{x=\infty\}, \\
 &\Gamma_c = \partial\Omega \cap \{\nu=0\}, && \Gamma_d = \partial\Omega \cap \{\nu=\infty\}.
\end{align*}
From the $\log$-transformed Heston PDE \eqref{eq: Heston transformed}, we define
\begin{align*}
A = \frac{1}{2} \nu \begin{pmatrix} \sn^2 & \sn \rho \\ \sn \rho & 1  \end{pmatrix}, \qquad 
b=\begin{pmatrix} \kn(\mr - \nu) -  \frac{1}{2} \sn^2 \\r-q-\frac{\nu}{2} - \frac{1}{2} \sn \rho \end{pmatrix}
\end{align*}
and rewrite it to its divergent form
\begin{equation*}
    V_\tau - \nabla \cdot A \nabla V - b \cdot \nabla V + rV = 0.
\end{equation*}
The \textit{constraint operator} $e$ in the Lagrangian is implicitly defined by  
\begin{equation} \label{eq: Lagrangian constraint}
\begin{aligned}
    \langle e(V,\xi) , \psi \rangle &=\int_0^T \int_{\tOmega} V \Bigl[- \varphi_\tau - \nabla \cdot A^\top \nabla \varphi + b\cdot \nabla \varphi + (r+ \nabla\cdot b) \varphi \Bigr] \,dz   \\
&\qquad+ \int_{\partial \Omega} [ (A^\top\nabla \varphi) \cdot \vec n  - (b\cdot \vec n) \varphi]V \,ds
      - \int_{\partial \Omega} (A \nabla V ) \cdot \vec n \varphi \,ds d\tau
\\ 
    &\qquad + \int_0^T \int_{\Gamma_a} \bigl[ V - \exp(-r\tau) \bigr] \varphi^a \,ds d\tau
       + \int_0^T \int_{\Gamma_b} V  \varphi^b\, ds d\tau \\
    &\qquad + \int_0^T \int_{\Gamma_d} \bigl[  V - \exp(-r\tau)\bigr] \varphi^d\, ds d\tau,
\end{aligned}
\end{equation} 
where $\tilde{\Omega}=\Omega \cap \Gamma_c$, since no boundary condition is required at $\Gamma_c$.

Following the approach of the authors \cite{clevenhaus24gradient}, we obtain the adjoint equation
\begin{multline}\label{eq:HestonAD_noA}
     \varphi_{\tau} + \frac{1}{2} \nu \sn^2 \varphi_{\nu \nu}+ \nu \sn\rho \varphi_{x \nu} + \frac{1}{2} \nu \varphi_{xx}+ (\sigma_\nu^2-\kappa_\nu(\mr - \nu) ) \varphi_{\nu} \\
     + (q-r+\frac{\nu}{2} +\sigma_\nu \rho) \varphi_{x} + (\kn -r) \varphi = -(V-\Vd) \quad \text{on }\Omega
\end{multline}
with terminal condition $\varphi(T)=0$ and $\varphi=0$ on the boundaries $\Gamma_a$, $\Gamma_b$ and $\Gamma_d$ and the outflow boundary at $\nu=0$. 
For a more detailed derivation of the adjoint formulation, see \cite{clevenhaus24gradient}.

We compute the optimality condition by setting $d_{\xi} L(V,\xi,\psi) = 0$. 
Since the boundaries $\Gamma_a$, $\Gamma_b$ and $\Gamma_d$ are zero, we focus on $\tOmega$. 
In the following, the derivatives with respect to the different parameters are given separately
\begin{align*}
 d_{\sn} \langle e(V,\xi),\psi \rangle  &= \int_0^T \int_{\tilde{\Omega}}  V \Bigl[ - \sn \nu \varphi_{\nu \nu} - 2 \sn \varphi_{\nu} -\rho \varphi_x - \rho \nu \varphi_{x \nu}\Bigr]\, dz\, d\tau,\\   
    d_{\rho} \langle e(V,\xi),\psi \rangle  
    &=   \int_0^T \int_{\tilde{\Omega}} V \Bigl[ - \sn \varphi_x - \sn \nu \varphi_{x \nu}\Bigr]\,dz\, d\tau, \\
    d_{\kn} \langle e(V,\xi),\psi \rangle 
    &= \int_0^T \int_{\tilde{\Omega}} V \Bigl[ (\mr-\nu)\varphi_{\nu}-\varphi \Bigr] \,dz\, d\tau, \\
     d_{\mr} \langle e(V,\xi),\psi \rangle 
     &=  \int_0^T \int_{\tilde{\Omega}}  \kn V \varphi_{\nu} \, dz\, d\tau.
\end{align*}
Note that $d_{\xi} L(V,\xi,\psi)[h_{\xi}] = 0$ must hold for arbitrary directions $h_{\xi}$. 
Therefore, we can read off the gradient from the above expressions.
We solve the state equation for the Heston model with a constant control variable $\xi$ \eqref{eq: Heston transformed} for a given initial parameter set $\xi_0$,
With the state solution at hand, we can compute the corresponding adjoint equation \eqref{eq:HestonAD_noA} backwards in time. The state and adjoint solutions allow us to compute the gradient for the descent step.

We ensure that $\kappa_\nu$, $\mu_\nu$, $\sigma_\nu$, and $\rho$ are within their parameter bounds, and that the Feller condition is satisfied by using the \textit{projected Armijo rule} \cite{Troeltzsch}.
In the projected Armijo rule, we choose the maximum $\sigma_k\in\{1,1/2,1/4,\ldots\}$ for which
\begin{equation*}
 f\bigl(\mathcal{P}(\xi_k - \sigma_k \nabla f(\xi_k))\bigr) - f(\xi_k) \le -\frac{\gamma}{\sigma_k}\|\mathcal{P}(\xi_k-\sigma_k \nabla f(\xi_k))-\xi_k\|_2^2.
\end{equation*}
Here, $\gamma\in(0,1)$ is a numeric constant, we will use the default choice with $\gamma=10^{-4}$.

%%%%%%%%%%%%%%%%%%%%%%%%%%%%%%%%%%%%%%%%%%%%%%%%%%% SECTION 4
\section{Numerical Scheme}\label{sec: Numerical scheme} 
Since this is a proof-of-concept, simple and well-known spatial and temporal discretization methods are used to illustrate our approach for discretizing both the SDE and the PDE.
%%%%%%%%%%%%%%%%%%%%%%%%%%%%%%% Section 4.1.
\subsection{The Fine Model: Heston SDE}
For the fine model, we consider a Monte Carlo simulation with variance reduction. 
We consider the $\log$-transformed SDE 
\begin{equation}
    \begin{cases}
 dx_t = (r-q-\half \nu_t)\,dt + \sqrt{\nu_t} \,d\bm{t}{S},&\quad S_0>0,\\
 d\nu_t = \kn(\mr-\nu_t)\,dt+\sn \sqrt{\nu_t}\,d\bm{t}{\nu},&\quad \nu_0>0, \\
d\bm{t}{x}\,d\bm{t}{\nu}=\rho\,dt
\end{cases}.
\end{equation}
and use the well-known Euler-Maruyama scheme to discretize the SDE.
Therefore we discretize the time uniformly $t_k=k \Delta t$, where $\Delta_t=T/N_t$ and $k=0,\ldots,N_t$.
Given two independent normal distributed random variables for each time step $\omega_{1,k}$ and $\omega_{2,k}$, we generate
\begin{equation}
\begin{aligned}
    \omega_{x,k}&=\omega_{1,k}, \\
    \omega_{\nu,k}&=\rho \omega_{x,k} + \sqrt{1-\rho^2}\omega_{2,k},
\end{aligned}    
\end{equation}
and obtain the discretized scheme 
\begin{equation}
    \begin{cases}
        x_{k+1}=x_k + (r-q-\half \nu_t)\,\Delta_t + \sqrt{\nu_t} \sqrt{\Delta t}\omega_{x,k}\\
        \nu_{t+1} = \nu_t + \kn(\mr-\nu_t)\,\Delta_t+\sn \sqrt{\nu_t}\sqrt{\Delta_t} \omega_{\nu,k}
    \end{cases}.
\end{equation}
To reduce the variance, we use antithetic variables \cite{Gunther.2010} to obtain the so-called antithetic path $x^{-}$. 
Since we use a uniform spacing for the time discretization, the discrete integral $A(0,T)$ of the payoff function is approximated by
\begin{equation}
   A(0,T)=\frac{1}{T}\int_0^T \exp(x_t) dt \approx \frac{1}{T} \frac{T}{N_t}\sum_{k=1}^{N_t}  \exp(x_k) = \frac{1}{N_t} \sum_{k=1}^{N_t}\exp(x_k)=\mean(x).
\end{equation}
If $N_p$ is the number of paths generated for the Heston model, we get the payoff for a path $p$
\begin{equation}
   \phi_p(x, x^-,T)= \frac{1}{2} \Bigl( \max\bigl(K-\mean(x),0\bigr) 
                            + \max\bigl(\bigl(K-\mean(x^-),0\bigr) \Bigr)
\end{equation}
and the fair price of the option is given by 
\begin{equation}
    \exp{(-rT)} \frac{1}{N_p}\sum_{i=1}^{N_p}\phi_{p_i}(x, x^-,T).
\end{equation}

%%%%%%%%%%%%%%%%%%%%%%%
\subsection{The Coarse Model: Heston PDE}
Before discretizing the Heston PDE, we perform a domain truncation and introduce a closure boundary condition at $\nu=0$ for the Heston and its adjoint formulation. 
Since we assume that the Feller condition holds, according to Fichera theory we have a pure outflow boundary at $\nu=0$ and no need for an analytical boundary condition, neither for the Heston model nor for its adjoint formulation.
However, we do need a closure condition at this boundary for the discretized PDE. 
As discussed in \cite{clevenhaus2023ECMI, Kutik.2015}, we follow Heston's approach and use the reduced hyperbolic formulation of the Heston PDE and its adjoint. 
At $\Gamma_c$, we obtain for the $\log$-transformed normalized PDE and similarly for the adjoint formulation 
\begin{equation}\
    \varphi_{\tau} + (\sigma_\nu^2-\kappa_\nu\mr) \varphi_{\nu} + (q-r+\sigma_\nu \rho) \varphi_x  + (\kn -r) \varphi = -(V-\Vd).
\end{equation}
We obtain a rectangular grid by truncating the domain and introducing  grid points. 
We consider uniform meshes in each direction and get $x_i=x_{\min}+i \Delta_x$ for $i=0,\ldots,N_x$ 
with $\Delta_x=(x_{\max}-x_{\min})/N_x$ 
and $\nu_j=j \Delta_\nu$ for $j=0,\ldots,N_\nu$ with $\Delta_\nu=\nu_{\max}/N_\nu$ for the spatial direction,
and $\tau_k=k\,\Delta_\tau$ for $k=0,\ldots,N_\tau$ with $\Delta_\tau=T/N_\tau$ for the time direction.

For the time discretization we use the Hundsdorfer-Verwer scheme \cite{Hundsdorfer.2002}, a well-known \textit{alternating direction implicit} (ADI) method that is able to handle mixed derivative terms.
Within the Hundsdorfer-Verwer scheme, $\theta$ is introduced as a classification measure similar to the $\theta$ method and is of second order for any choice of $\theta$.
Since the $\log$-transformed Heston PDE and its adjoint formulation are second-order PDEs, they are discretized similarly. 
Therefore, we introduce a general second-order PDE formulation
\begin{equation}
    u_\tau + a_{11} u_{\nu\nu} +2 a_{12} u_{x\nu} + a_{22} u_{xx} + b_1 u_{\nu} +b_2 u_x + c u=0.
\end{equation}

As a first step, we split the PDE operator into three operators
\begin{equation}   \label{eq: Operator Splitting ADI}
    \F(\tau) = \F_0(\tau) + \F_1(\tau) + \F_2(\tau).
\end{equation}
Each operator deals with a spatial direction, $F_0$ includes the mixed derivative terms, $F_1$ includes the derivatives in $x$-direction, and similarly $F_2$ includes the derivatives with respect to $\nu$. 
Since the term $cu$ is independent, it is split in half and added to both $F_1$ and $F_2$. 
We get the splitting
\begin{align*}
    \F_0(\tau) &= 2 a_{12} D_{x \nu}, \\
    \F_1(\tau) &= b_2 D_x + a_{22} D_{xx} - \half c I,\\
    \F_2(\tau) &= b_1 D_{\nu} +  a_{11} D_{\nu \nu} -\half c I,
\end{align*}
where $D_x$ describes the discretization matrix of the first derivative w.r.t.\ $x$ and correspondingly $D_{xx}$ of the second derivative w.r.t.\ $x$, $D_\nu$ and $D_{\nu\nu}$ of the first and second derivatives w.r.t.\ $\nu$, $I$ denotes the identity matrix. 
We use central finite difference stencils to derive the corresponding matrices.

Let $u^{i,j}_k\approx u(x_i,\nu_j,\tau_k)$ and simplify $u_k\approx u(x,\nu,\tau_k)$.
In each time step, we have to solve the following system of equations
\begin{equation}
    \begin{cases}
    Y_0=u_{k}+\Dtau \F(\tau_{k}) u_{k}, \\
    Y_1=Y_{0}+\theta \Dtau \bigl(\F_1(\tau_{k+1})Y_1 -\F_1(\tau_{k}) u_{k} \bigr), \\
    Y_2=Y_{x}+\theta \Dtau \bigl(\F_2(\tau_{k+1},Y_2)-\F_2(\tau_{k}) u_{k} \bigr), \\
    \tY_0=Y_0+\half \Dtau \bigl(\F(\tau_{k+1})Y_2-\F(\tau_{k}) u_{k}\bigr),\\
    \tY_1=\tY_{0}+\theta \Dtau \bigl(\F_1(\tau_{k+1})\tY_1-\F_1(\tau_{k})u_{k}\bigr), \\
    \tY_2=\tY_{x}+\theta \Dtau \bigl(\F_2(\tau_{k+1})\tY_2-\F_2(\tau_{k})u_{k}\bigr), \\
    u_{k+1}=\tY_2,
    \end{cases}
\end{equation}
where we choose $\theta=3/4$ and implement $u_k$ as a matrix \cite{Teng.2019}. 
The discrete boundary conditions in $x$-direction for a sufficiently small $x_{\min}$ and a sufficiently large $x_{\max}$ are set to $V(x_{\min},\nu,\tau)=\exp(-r\tau)$ and similarly $V(x_{\max},\nu,\tau)=0$. 
For the variance at $\nu=\nu_{\min}$ and $\nu=\nu_{\max}$ we use the ghost cell approach. 
Within this approach, additional grid points are computed at $\nu_0=\nu_{\min}-\Delta_{\nu}$ and $\nu_{N_{\nu}+2}=\nu_{\max}+\Delta_{\nu}$ via zero-order extrapolation. 
For a detailed discussion, see \cite{clevenhaus2023ECMI, Kutik.2015}.
The integrals appearing in the gradient are computed by the trapezoidal rule.

%%%%%%%%%%%%%%%%%%%%%%%%%%%%%%%%%%%%%%%%%%% Numerical Results
\section{Numerical Results}\label{sec: Numerical Results}
From the put options on the Nikkei~300 index on Dec. 31,
2012, we get one $S_0$ and five different sets, each with the following parameters $r$, $q$, and $K$.
Thus, we specify $M1, M2, M3, M4$ and $M5$ for the different market data sets.
Since the maturity must be small to use the European option price as a proxy for the Asian option price, we focus on $T=0.25$.
For the spatial discretization we set 
\begin{equation}
\begin{split}
       x_{\max}=\log(1.2*S_0),\quad N_x=120,\quad \Delta_x=\frac{x_{\max}}{N_x} \quad\text{and}\quad x_{\min}=\Delta_x, \\
    \nu_{\max}=1, \quad N_{\nu}=100, \quad \Delta_x=\frac{\nu_{\max}}{N_{\nu}} \quad\text{and} 
    \quad\nu_{\min}=\Delta_{\nu}. 
\end{split}
\end{equation}
and select $\nu_0=0.05$. 
As a result of this discretization, the strike price (and thus the kink in the payoff function) is approximated at a grid point. 
Therefore, we smooth the initial condition using the operator from Kreiss et al.~\cite{Kreiss}
The time discretization uses
\begin{equation}
    \tau_{\max}=T,\quad 
    N_{\tau}=170,\quad 
    \Delta_{\tau}=\frac{T}{N_{\tau}} \quad
    \text{and}\quad \tau_{\min}=0. \\
\end{equation}
We consider six different initial guesses for the algorithm, 
see Table~\ref{tab:init_guess}; each parameter set satisfies the Feller condition.  
\begin{table}[]
    \centering
    \begin{tabular}{|c|c|c|c|c|} \hline 
       Guess  & $\kn$ & $\mr$&$\sn$&$\rho$ \\ \hline
        1 &3.0& 0.3&0.1&-0.2  \\ \hline
        2 &5.0& 0.6&0.2&-0.3 \\ \hline
        3 &4.5&0.8&0.5&-0.15  \\ \hline
        4 &2.0&0.4&0.45&-0.2  \\ \hline
        5 &4.0& 0.5&0.15&-0.35 \\ \hline
        6 &3.5& 0.35&0.5&-0.5  \\ \hline
    \end{tabular}
    \caption{Different initial guess sets for the initial coarse model calibration.}
    \label{tab:init_guess}
\end{table}

First, we focus on the gradient descent algorithm. For the algorithm, we set the iteration maximum to 51 and the terminal condition $J_c(V,c(\xi_{\mathrm{c}}))<10^{-3}$. 
Since we are using a gradient-based algorithm, we can only expect to converge to a local minimum, to evaluate the descent over the iterations we use the relative reduction of the cost functional 
\begin{equation} \label{eq: rel_red}
    r(\xi_{\mathrm{c}}^0)= 100 \cdot \left(1-\frac{J(V(\xi_{\rm opt});\Vd)}{J(V(\xi_{\mathrm{c}}^0);\Vd)}\right),
\end{equation}
where $\xi_{\mathrm{c}}^0$ is the initial guess from Table~\ref{tab:init_guess}.
%because the algorithm is sensitive to the initial guess.

\begin{figure}
    \centering
    \includegraphics[width=\linewidth]{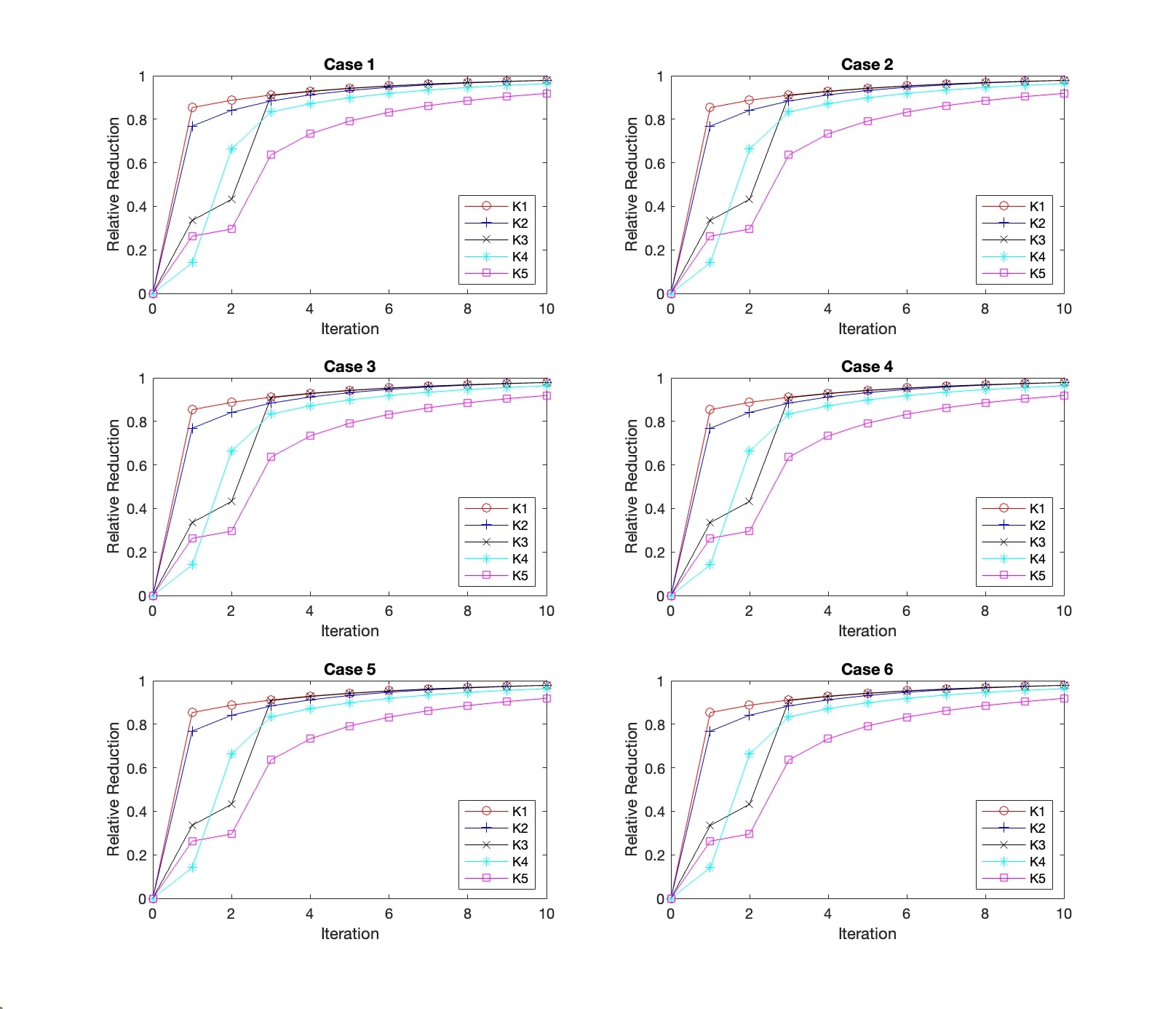}
    \caption{Cost functional reduction for the first 10 Iterations for the initial calibration for the coarse model.}
    \label{fig:coarse}
\end{figure}

Table~\ref{tab:coarse_red} shows the cost functional reduction of the last and/or optimal cost function value.
The bold values indicate instances where the maximum number of iterations was reached.
We observe that the gradient descent algorithm calibrates the parameter $\xi_{\mathrm{c}}$ almost perfectly, even if we can't guarantee to find the global minimum. 
In Table~\ref{tab:coarse_cost}, which shows the optimal value of the cost functional, 
we observe that for most of the test cases we reach the terminal condition of the gradient algorithm.
The cases where the condition is not reached at the iteration maximum are shown in bold.
The cases with the highest cost functional values correspond to the smallest cost function reduction.
To illustrate the cost functional reduction per iteration, Figure~\ref{fig:coarse} shows the value for the first 10 iterations. 
We observe that the maximum number of iterations can be significantly reduced depending on the desired accuracy.

\begin{table}[]
    \centering
    \begin{tabular}{|c|c|c|c|c|c|c|} \hline
      M \textbackslash Guess  &1 &  2 & 3 & 4 & 5 & 6  \\ \hline
     1& 99.97 &  99.91 & 99.98 & 99.66 & 99.99 & \textbf{99.95} \\  \hline
     2& 99.97 &  99.92 & 99.98 & 99.81 & 99.98 & 99.96  \\ \hline
    3&  99.97&  \textbf{99.51} & \textbf{99.69} & 99.50 & 99.98 & \textbf{99.93}  \\  \hline
     4&   99.97 &  \textbf{99.75} & \textbf{99.97} & 99.72 & 99.96 & \textbf{99.87}  \\ \hline
      5&   \textbf{99.96} &  99.92 & 99.98 & 99.92 & \textbf{96.15} & \textbf{99.49}  \\ \hline
   \end{tabular}
    \caption{Cost functional reduction for the initial calibration of the coarse model.}
    \label{tab:coarse_red}
\end{table}
\begin{table}[]
    \centering
    \begin{tabular}{|c|c|c|c|c|c|c|} \hline
      M \textbackslash Guess  &1 &  2 & 3 & 4 & 5 & 6  \\ \hline
     1& 0.888 &  0.861 & 0.578 & 0.808 & 0.289 & \textbf{1.058}  \\  \hline
     2& 0.928 &  0.845 & 0.680 & 0.603 & 0.401 & 0.892  \\ \hline
    3&  0.933 &  \textbf{4.638} & \textbf{9.914} & 0.516 & 0.461 & \textbf{1.702}  \\  \hline
     4&   0.906 &  \textbf{1.639} & \textbf{1.248} & 0.678 & 0.707 & \textbf{3.542}  \\ \hline
      5&   \textbf{1.031} &  0.928 & 0.988 & 0.689 & \textbf{56.160} & \textbf{17.323}  \\  \hline
   \end{tabular}
    \caption{Cost functional value for the optimal calibration of the coarse model, scaled by $10^3$.}
    \label{tab:coarse_cost}
\end{table}

These results show that the presented gradient descent algorithm is a viable choice for calibration. For a more detailed analysis of the gradient descent algorithm, see \cite{clevenhaus24gradient}. 
Within the space mapping, we limit ourselves to a maximum of four iterations, since the evaluation of the fine model is expensive.
As a calibration measure we again use the relative cost functional reduction, the results are presented in Table~\ref{tab:sm_red}. 
The results show that the space mapping approach has only slightly lower reduction rates, except for K5. 
One can improve this value by choosing a smaller $\nu_0$, as can be observed in K5a, where we choose $\nu_0=0.03$. 
Note that the choice of $\nu_0$ is limited by the grid structure.
In addition, the calibration to $\nu_0$ can be included in the gradient descent algorithm. 

Table~\ref{tab:sm_cost} shows the overall cost functional reduction relative to the initial guess and Table~\ref{tab:sm_red} shows the optimal cost functional value.
Similar to the gradient algorithm, we observe that when the cost functional reduction is small, the cost functional value is larger. 

\begin{table}[]
    \centering
    \begin{tabular}{|c|c|c|c|c|c|c|} \hline
      M \textbackslash Guess &1 &  2 & 3 & 4 & 5 & 6  \\ \hline
     1& 98.57 &  98.75 & 99.07 & 96.73 & 99.52 & 99.89 \\  \hline
     2& 99.66 &  99.21 & 99.96 & 72.41 & 98.93 & 99.22  \\ \hline
    3&  98.01&  99.64 & 99.98 & 98.33 & 91.67 & 99.39  \\  \hline
     4&   92.42 &  99.20 & 93.42 & 98.44 & 98.05 & 99.74  \\ \hline
      5&   28.22 &  53.71 & 64.48 & 52.74 & 42.48 & 71.84  \\  \hline
      5a& 44.01 & 46.58& 43.52& 61.80& 73.26& 79.14\\ \hline
   \end{tabular}
    \caption{Cost functional reduction for the calibration of the space mapping.}
    \label{tab:sm_red}
\end{table}

\begin{figure}
    \centering
    \includegraphics[width=\linewidth]{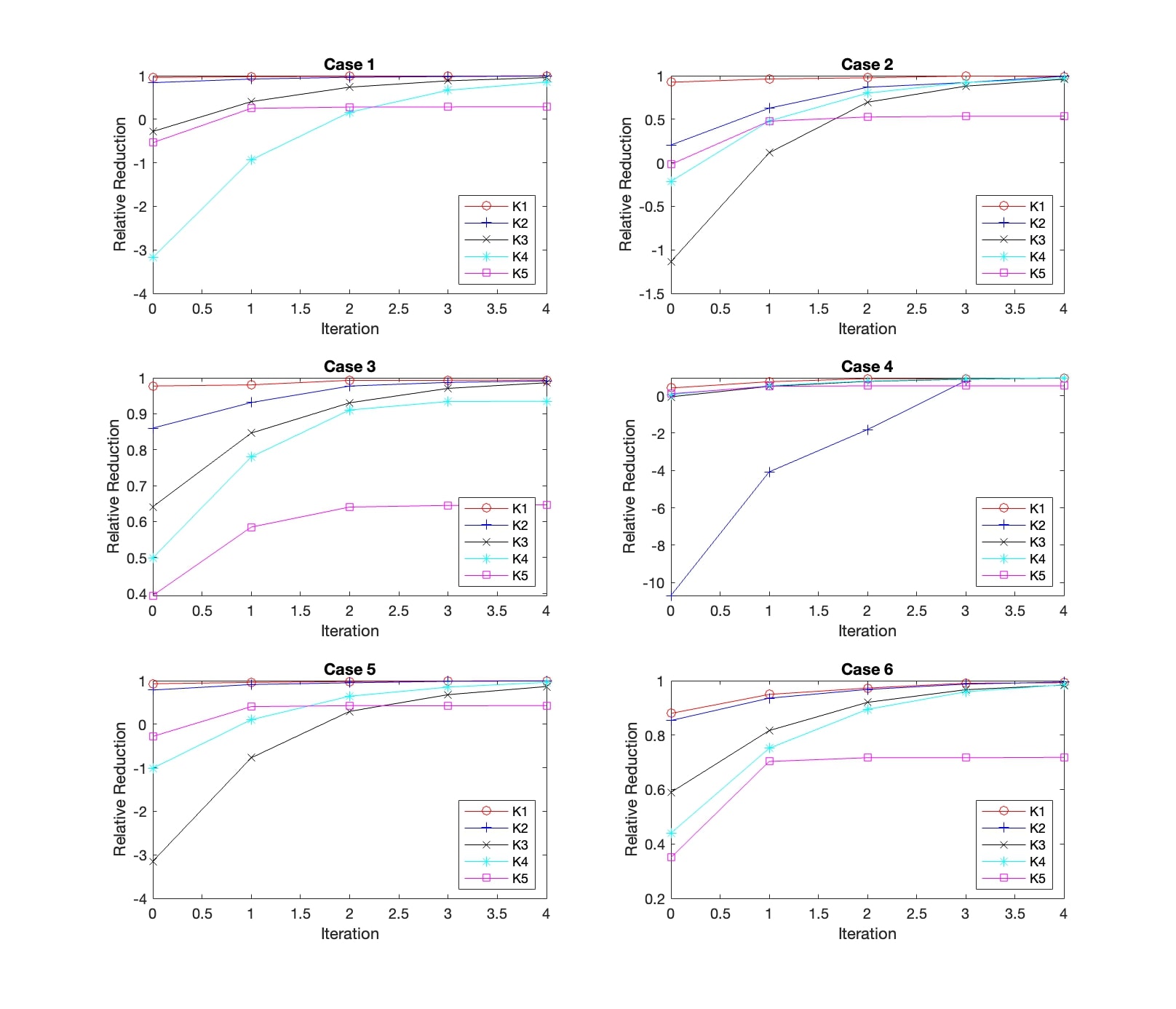}
    \caption{Cost functional reduction within the space mapping algorithm for the fine model adjusted to $\xi_{\rm{init}}$.}
    \label{fig:fine_guess}
\end{figure}

\begin{figure}
    \centering
    \includegraphics[width=\linewidth]{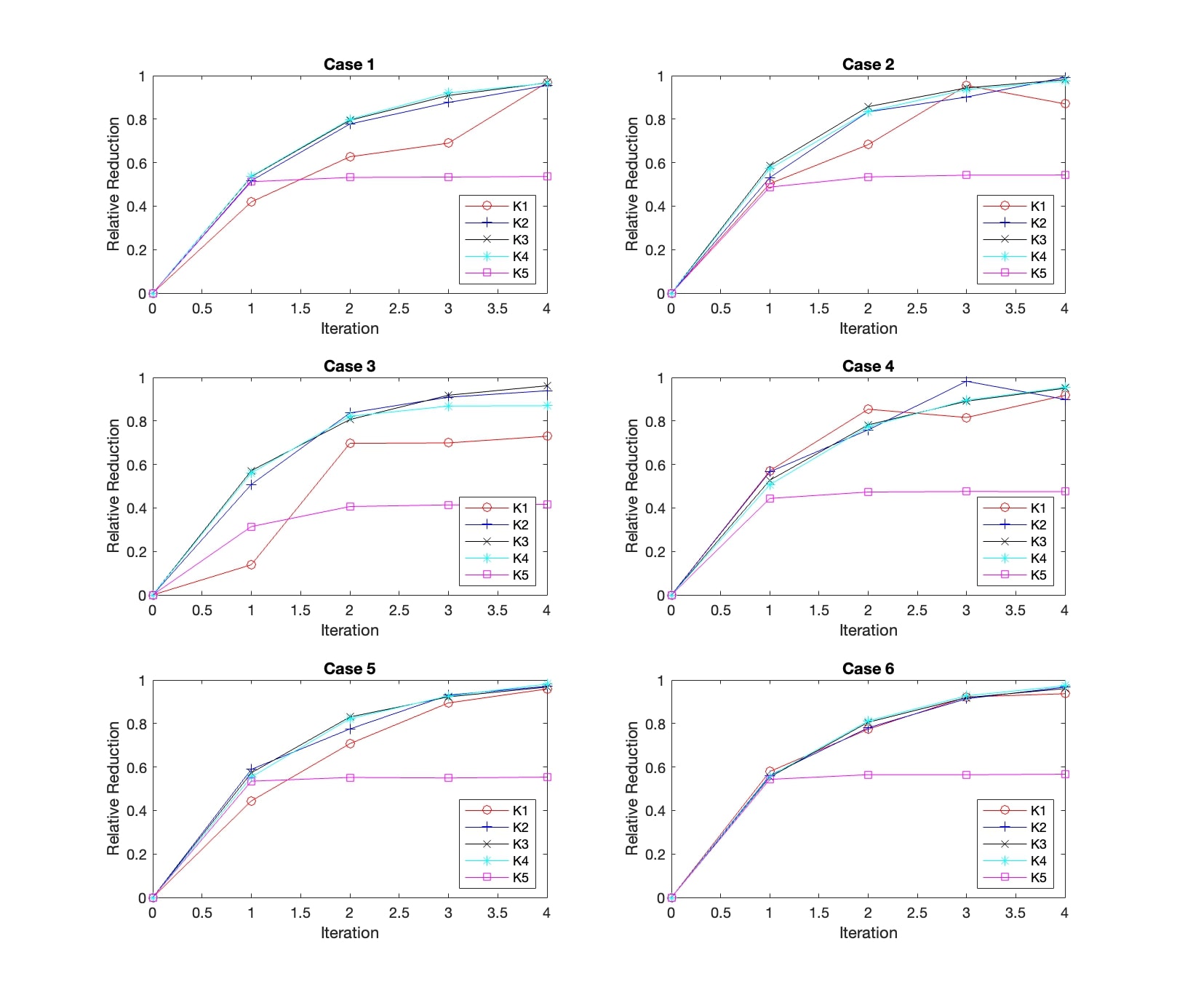}
    \caption{Cost functional reduction within the space mapping algorithm for the fine model adjusted to $\xi_f^0$.}
    \label{fig:fine_init}
\end{figure}

Since $J(V_c(\xi_c^*);\Vd)\approx J_f(V_f(\xi_f^0);\Vd)$ can be significantly worse than $J(V_f(\xi_{\mathrm{c}}^0);\Vd)$, we present two figures. One figure shows the cost functional reduction per iteration, w.r.t.\ the initial guess $\xi_{\mathrm{c}}^0$, see Figure~\ref{fig:fine_guess}, and the other one w.r.t.\ the optimal calibration parameters resulting from the initial coarse model calibration $\xi_c^*=\xi_f^0$, see \ref{fig:fine_init}. 
The Figures~\ref{fig:fine_guess} and \ref{fig:fine_init} and Table~\ref{tab:grad_red_A} show that even if $\xi_f^0$ results in a higher cost functional value for the space mapping at iteration~0, 
the space mapping reduces the cost functional significantly, e.g.\ Guess~1 with $M4$, Guess~4 with $M2$ as well as Guess~5 with $M3$.
\begin{table}[]
    \centering
    \begin{tabular}{|c|c|c|c|c|c|c|} \hline
      M \textbackslash Guess  &1 &  2 & 3 & 4 & 5 & 6  \\ \hline
     1& 0.021 &  0.006 & 0.012 & 0.007 & 0.005 & 0.001  \\  \hline
     2& 0.004 &  0.002 & 0.001 & 0.003 & 0.009 & 0.009  \\ \hline
    3&  0.012 &  0.001 & 0.022 & 0.011 & 0.015 & 0.0110  \\  \hline
     4&   0.027 &  0.010 & 0.2010 & 0.024 & 0.014 & 0.007  \\ \hline
      5&   1.159 &  1.124 & 1.480 & 1.290 & 1.110 & 1.046  \\  \hline
      5a& 0.574& 0.641& 1.001& 0.699& 0.593& 0.492\\ \hline
   \end{tabular}
    \caption{Cost functional value for the optimal calibration of the space mapping.}
    \label{tab:sm_cost}
\end{table}
\begin{table}[]
    \centering
    \begin{tabular}{|c|c|c|c|c|c|c|} \hline
      M  \textbackslash Guess  &1 &  2 & 3 & 4 & 5 & 6  \\ \hline
     1& 96.12 & 92.53 & 97.73& 41.27 & 93.00 & 87.99  \\  \hline
     2& 84.12 & 20.39 & 86 01& -1069.14 & 78.50 & 85.35  \\ \hline
     3& -28.13 & -113.55 & 64.14& -5.97 & -314.86 & 58.95  \\  \hline
     4& -317.71 & -21.38 & 49.87 & 5.85 & -100.91 & 43.96  \\ \hline
     5&  -53.76&  -1.80&  39.41 & 10.54 & -28.04  & 35.00\\ \hline
   \end{tabular}
    \caption{Cost functional reduction for the initial guess and the calibration of the coarse model.}
    \label{tab:grad_red_A}
\end{table}

Since the purpose of this article is to introduce space mapping to computational finance research as a proof-of-concept, no performance analysis is provided. 
However, by using faster or more accurate techniques for solving the Heston PDE and its adjoint equations, as well as for Heston calibration, the algorithm can be optimized.

%%%%%%%%%%%%%%%%%%%%%%%%%%%%%%%%%%%%%%%%%%% Conclusion
\section{Conclusion and Outlook}\label{sec: Conclusion and Outlook}
The numerical results show the feasibility of the space mapping approach as a new calibration method in financial research.
The numerical results are remarkable even in this proof-of-concept article and can be improved in several ways, e.g., by including time-dependent parameters, adding a calibration for $\nu_0$, and using improved computational methods for the coarse model optimization and for the fine model solver.
In particular, the gradient descent algorithm can be easily adapted to account for time-dependent parameters to improve space mapping for pricing Asian options since they are time-dependent, see \cite{clevenhaus24gradient}.
The space mapping approach can be applied to various other hierarchical problems in finance, e.g.\ model, temporal and spatial as well as option hierarchies. 
Therefore, this article is a first step towards integrating space mapping into financial applications.

%%%%%%%%%%%%%%%%%%%%%%%%%%%%%%%%%%%%%%%%%%% Refs
\bibliographystyle{plain}
%\bibliography{SpaceMapping/SM}

\end{document}